\documentclass{amsart}
\newif\ifpdf \ifx\pdfoutput\undefined \pdffalse \else \pdfoutput=1 \pdftrue \fi
\ifpdf
   \usepackage[pdftex]{hyperref}
   \usepackage{color,graphicx}
   \DeclareGraphicsExtensions{.pdf,.jpg}
\else
   \usepackage{color,graphicx}  % Graphics extensions in graphics are eps by default
\fi
\ifpdf

\else

\fi
\begin{document}
\bibliographystyle{plain}
\newtheorem{theorem}{Theorem}
\newtheorem{lemma}{Lemma}
\newtheorem{remark}{Remark}
\newtheorem{conj}{Conjecture}
\newtheorem{prop}{Proposition}
\newtheorem{cor}{Corollary}
\newtheorem{df}{Definition}
\date{}
\newcommand{\complaint}[1]{ \framebox[1.2\textwidth]{\parbox{1.1\textwidth}{#1}} \\ }
%%%%%%%%%%%%%%%%%%%%%%%%%%%%%%%%%%%%%%%%%%%%%%%%%%%%%%%%%%%%%%%%%%%%%%%
\title{Combinatorial pseudo-Triangulations}
\author[Orden]{
     David Orden}
\thanks{Work of D. Orden and F. Santos was partially supported by grant
BFM2001--1153 of the Spanish  Direcci\'on General de Ense\~nanza
Superior e Investigaci\'on Cient\'{\i}fica}
\address{
Departamento de Matem\'aticas, Universidad de Alcal\'a. Facultad
de Ciencias, Apdo. de Correos 20, E-28871 Alcal\'a de Henares
(Madrid), SPAIN} \email{david.orden@uah.es}
\author[Santos]{
     Francisco Santos}
%\thanks{Work of D. Orden and F. Santos was partially supported by grant
%BFM2001--1153 of the Spanish  Direcci\'on General de Ense\~nanza
%Superior e Investigaci\'on Cient\'{\i}fica}
\address{
Departamento de Matem\'aticas, Estad\'{\i}stica y Computaci\'on,
Universidad de Cantabria. Facultad de Ciencias, Avenida de los
Castros s/n, E-39005 Santander, SPAIN } \email{santosf@unican.es}
\author[Servatius]{
     Brigitte Servatius \and
     Herman Servatius}
\thanks{Work of B. Servatius and H.\ Servatius was partially supported by
NSF grant 0203224.}
\address{Dept. of Math. Sci., Worcester Polytechnic Institute, Worcester MA 01609}
\email{bservat@wpi.edu, hservat@wpi.edu}
%%%%%%%%%%%%%%%%%%%%%%%%%%%%%%%%%%%%%%%%%%%%%%%%%%%%%%%%%%%%%%%%%%%%%%%

\begin{abstract}
   We prove that a planar graph is generically rigid in the plane if and only
   if it can be embedded as a pseudo-triangulation. This generalizes
  the main result of~\cite{horsssssw} which treats the minimally generically
  rigid case.

  The proof uses the concept of combinatorial pseudo-triangulation, CPT, in the plane
   and has two main steps: showing that a certain ``generalized Laman
property'' is a necessary and sufficient condition for a CPT to be ``stretchable'', and
showing that all generically rigid plane graphs admit a CPT assignment with that property.

Additionally, we propose the study of combinatorial pseudo-triangulations on
closed surfaces.
\end{abstract}

\maketitle

\section{Introduction}
\label{sec:intro}
The study of
pseudo-triangulations
%(also called geodesic triangulations~\cite{SomethingShouldBeCited}, in the case of polygons)
in the plane was initiated recently \cite{cegges,pv}, but it is rapidly
becoming a standard  topic in Computational Geometry.
This paper continues the program established in~\cite{horsssssw,os-pncgp-02,rss,streinu},
showing that {\em pseudo-triangulations provide a
missing link between planarity and rigidity in geometric graph theory}.
Our main result is:

\begin{theorem}\label{thm:main}
For a plane graph $G$, the following conditions are equivalent:
\begin{enumerate}
\item $G$ is generically rigid in the plane.
\item $G$ can be stretched to become a pseudo-triangulation
of its vertex set with the given topological embedding.
\end{enumerate}
\end{theorem}

Recall that a \emph{plane graph}
is a graph together with a given (non-crossing) topological embedding in the plane.
It is interesting to observe that property (1) is a property of the underlying (abstract) graph,
while property (2) in principle depends on the embedding. That is, our result in
particular implies that pseudo-triangulation stretchability of a plane graph is independent
of the embedding. Another consequence is that:

\begin{cor}\label{coro:main}
The class of planar and generically rigid (in the plane) graphs
coincides with the class of graphs of pseudo-triangulations.
\end{cor}

The implication from (2) to (1) in Theorem~\ref{thm:main} follows
from~\cite{os-pncgp-02}. Here we prove the implication from (1) to
(2) in two steps, via the concept of {\em combinatorial
pseudo-triangulation } introduced in Section~\ref{sec:cpts}.
Section~\ref{sec:Laman-labelling} proves that generically rigid
plane graphs can be turned into combinatorial
pseudo-triangulations with a certain ``generalized Laman''
property, and Section~\ref{sec:CPT-stretch} proves that
combinatorial pseudo-triangulations with this property can be
stretched. In Section~\ref{sec:surfaces} we initiate the study of
combinatorial pseudo-triangulations on closed surfaces, which we
think is an interesting topic to be developed further.

%%%%%%%%%%%%%%%%%%%%%%%%%%%%%%%%%%%%%%%%%%%%%%%%%%%%%%%%%%%%
\section{Pseudo-triangulations and rigidity}
\label{sec:rigidity}

Let $A$ be a finite point set in the Euclidean plane, in general position.
A {\em pseudo-triangle}
in the plane is a simple polygon with exactly three convex angles. A pseudo-triangulation of $A$
is a geometric (i.e., with straight edges) non-crossing graph with vertex set $A$,
containing the convex hull edges of $A$ and in which
every bounded region is a pseudo-triangle. A vertex $v$ in a geometric graph $G$ is called
{\em pointed} if all the edges
of $G$ lie in a half-plane supported at $v$ or, equivalently, if one of the angles incident to $v$
is greater than $180^\circ$. A pseudo-triangulation is called {\em pointed} if all its vertices
are pointed. The following numerical result has been stated several times in different forms,
and is crucial to some of the nice properties of pseudo-triangulations:

\begin{lemma}
\label{lemma:pt-count}
Let $G$ be a non-crossing straight-line embedding of a connected graph in the plane.
Let $e$, $x$ and $y$ denote the numbers of edges, non-pointed vertices and pointed vertices
in $G$. Then, $e\le 3x+2y-3$, with equality if and only if the embedding is a pseudo-triangulation.
\end{lemma}

\begin{proof}
Let $f$ be the number of bounded faces of the embedding. By
Euler's formula, $x+y+f=e+1$. We now double-count the number of
``big'' and ``small'' angles in the embedding (that is, angles
bigger and smaller than 180 degrees, respectively). The total
number of angles equals $2e$. The number of big angles equals $y$,
and the number of small angles is at least $3f$ (every bounded
face has at least three corners) with equality if and only if the
embedding is a pseudo-triangulation. These equations give the
statement.
\end{proof}

Recall that a graph is {\em generically rigid in the plane}, see~\cite{gss} or~\cite{walterS}, if any generically embedded bar
and joint framework corresponding to the graph has no non-trivial infinitesimal motions.  Generic rigidity is a
property of the graph, and not of any particular embedding. In fact, edge-minimal generically rigid graphs on a given number
$n$ of vertices are characterized by {\em Laman's Condition}:  they have exactly $2n-3$
edges and every subset of $k$ vertices spans a subgraph with at most $2k-3$ edges, see~\cite{Laman}.
Generically rigid graphs with $|E| = 2n-3$ are also known as {\em Laman graphs}.

The connection between rigidity and pseudo-triangulations was first pointed out in
Streinu's
seminal paper~\cite{streinu} where it is proved that the graphs of pointed pseudo-triangulations
are minimally generically rigid graphs, i.e.\  Laman graphs.
In~\cite{horsssssw} it was shown that a graph $G$ has a realization as a pointed pseudo-triangulation
in the plane if and only if the graph is a planar Laman graph.
The following theorem in~\cite{os-pncgp-02} extends Streinu's result to non-minimally rigid graphs and
 relates the number of non-pointed vertices to the degree to which a planar rigid graph is
overbraced.
\begin{theorem}
\label{thm:ordensantos} Let $G$ be the graph of  a pseudo-triangulation of a planar point set
in general position.
Then:
\begin{enumerate}
\item $G$ is infinitesimally rigid, hence rigid and generically rigid.
\item Every subset of $x$ non-pointed plus $y$ pointed  vertices of $G$, with $x + y \geq 2$, spans
a subgraph with at most $3x + 2y -3$ edges.
\end{enumerate}
\end{theorem}

Property (2) will be crucial in our proof of Theorem~\ref{thm:main}.
 Observe, however,
that it is not a property of the graph $G$, but a property of the specific straight-line embedding of $G$.

Another remarkable connection between planarity, rigidity, and
pseudo-triangula\-tions concerns  planar rigidity circuits. These
are redundantly rigid graphs such that the removal of any edge
leaves a minimally rigid graph.  They can, by our results, be
realized as pseudo-triangulations with exactly one non-pointed
vertex.  Rigidity circuits, or Laman circuits, have the nice
property that the number of faces equals the number of vertices
and that their geometric dual (which exists and is unique) is also
a Laman circuit, see~\cite{ChristopherServatius}.  In~\cite{osssw}
we show, using techniques developed
here and Maxwell's classical theory of reciprocal
diagrams~\cite{Maxwell}, that if $C$ is a planar Laman circuit, then $C$ and
its geometric dual $C^*$  can be realized as
pseudo-triangulations with the same directions for corresponding
edges.

%%%%%%%%%%%%%%%%%%%%%%%%%%%%%%%%%%%%%%%%%%%%%%%%%%%%
\section{Combinatorial pseudo-triangulations in the plane}
\label{sec:cpts}

We now consider  a combinatorial analog of pseudo-triangulations.
Let $G$ be a plane graph.
We call {\em angles} of $G$ the pairs of consecutive edges in the
vertex rotations corresponding to the embedding. Equivalently, an angle
is a vertex-face incidence.
By a {\em labelling of angles} of $G$ we mean an assignment of
``big'' or ``small'' to every angle of $G$. Such a labelling is called
a {\em combinatorial pseudo-triangulation labelling} (or {\em
CPT-labelling}, for short) if every bounded face has exactly three
angles labelled ``small'', all the angles in the unbounded face are labelled
``big'', and no vertex is incident to more than one ``big'' angle.

\begin{figure}[htb]
\centering a.\includegraphics[scale=.60]{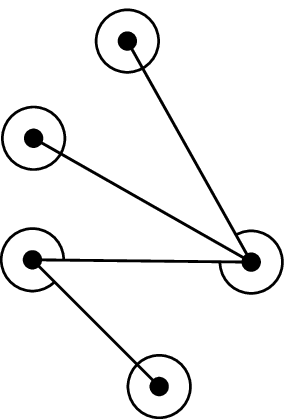} \qquad
b.\includegraphics[scale=.60]{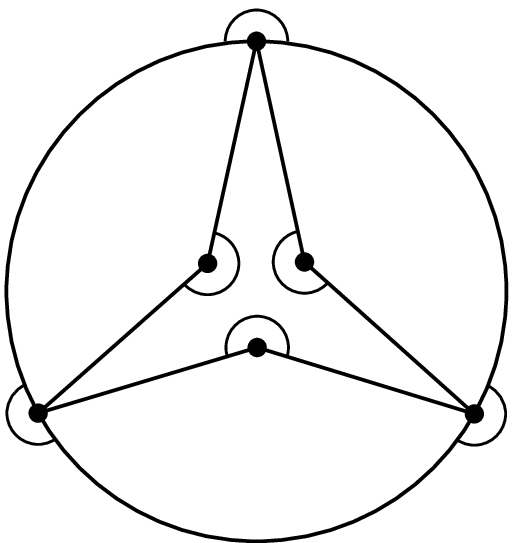} \qquad
c.\includegraphics[scale=.60]{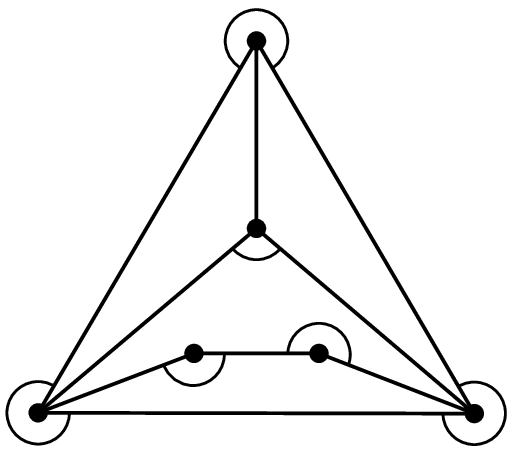} \caption{A tree with
unbounded face a pseudo-triangle (a), the Mercedes Graph as a CPT
(b), and a non-stretchable CPT (c).\label{FigMercedes}}
\end{figure}

The embedded graph $G$ together with a CPT-labelling of its angles
is called a {\em combinatorial pseudo-triangulation}, or CPT.
In figures we will indicate the large angles by an arc near the vertex between the edge pair.
Figure~\ref{FigMercedes} shows three graphs with large angles labelled. The one in the left
is not a CPT, because the exterior face has three small angles.
Figure~\ref{FigMercedes}b is a CPT, whose bounded faces are
three ``triangles'' and a hexagonal ``pseudo-triangle''.
If possible we shall draw large angles larger than $180^\circ$, small ones as angles smaller than $180^\circ$ and
edges as straight non-crossing segments, but it has to be observed
that this is sometimes not possible, since there are non-stretchable CPT's, such as the one in
Figure~\ref{FigMercedes}c.

Following the terminology of true pseudo-triangulations we say
that the interior faces of a CPT are {\em pseudo-triangles} with
the three small angles joined by three {\em pseudo-edges}. As in
the geometric case, a vertex is called {\em pointed} if there is a
big angle incident to it and the CPT is called {\em pointed} if
this happens at every vertex, see Figure~\ref{FigMercedes}b. The
following result and its proof, a straightforward counting
argument using Euler's formula, are completely analogous to the
geometric situation.

\begin{lemma}
\label{lemma:cpt-count}
Every combinatorial pseudo-triangulation in the plane with $x$ non-pointed  and $y$ pointed vertices
has exactly $3x+2y-3$ edges.
\end{lemma}

We recall here the main result of \cite{horsssssw}:

\begin{theorem} \label{thm:equivalences-lamangraphs}
     Given a plane graph $G$, the following conditions are equivalent:
     \begin{itemize}
     \item[(i)] $G$ is generically minimally rigid (isostatic),
     \item[(ii)] $G$ satisfies Laman's condition,
     \item[(iii)] $G$ can be labelled as a pointed CPT.
     \item[(iv)] $G$ can be stretched to a pointed pseudo-triangulation preserving the given topological
     embedding.
     \end{itemize}
\end{theorem}

Throughout the paper, we say that a combinatorial pseudo-triangulation $G$ has the
\emph{generalized Laman property} or is \emph{generalized Laman}
if every subset of $x$ non-pointed plus $y$ pointed  vertices,
with $x+y\geq 2$, induces a subgraph with at most $3x + 2y - 3$
edges. This property is inspired by~Theorem~\ref{thm:ordensantos} and it is
crucial to our proof; see Theorem~\ref{thm:main2}.
We call it the generalized Laman property because it
restricts to the Laman condition for the pointed case.

\begin{lemma}\label{lem:reformulated}
The generalized Laman property is equivalent to requiring that
every subset of $x'$ non-pointed plus $y'$ pointed vertices of
$G$, with $x'+y'\leq n-2$ be incident to at least $3x'+2y'$ edges.
\end{lemma}

\begin{proof}Using Lemma~\ref{lemma:cpt-count},
$x$ plus $y$ vertices satisfy the condition in the definition of
generalized Laman if and only if the cardinalities, $x'$ and $y'$,
of the complementary sets of vertices satisfy this reformulated
one.
\end{proof}

Lemma~\ref{lem:reformulated} implies that the generalized Laman
property forbids vertices of degree~1 and  that vertices of
degree~2 must be  pointed. Moreover, any edge cutset separating the
graph into two components, each containing more than a single
vertex, has cardinality at least 3.

%We will also frequently use
%the following lemma.
%
%
%\begin{lemma}\label{lem:extlamrig}
%   Let $G$ be a plane graph with a
%    CPT labelling satisfying the generalized Laman condition.
%   Then $G$ is generically rigid.
%\end{lemma}
%
%\begin{proof}
%     If there is a vertex of degree 2 then its deletion, together with
%     an easy relabelling, gives the result by induction.
%
%     Let  $x$ and $y$ be the number of non-pointed and pointed vertices
%     respectively. If $x=0$ then our statement follows from Theorem~\ref{thm:equivalences-lamangraphs}.
%     Otherwise, since $y>0$ (because all angles in the outer face are big) and since $G$ is connected,
%     there is an edge in the CPT with one pointed and one non-pointed
%     endpoint.  Removing such an edge, relabelling the merged small angles at the non-pointed
%     endpoint as big, and labelling the merged angles at the pointed endpoint big or small depending on
%     whether one of the angles to be merged is big or not, yields a graph with a generalized Laman
%     CPT.  Proceeding in this manner until all vertices are pointed produces a Laman graph.
%\end{proof}

The following is a more detailed formulation of our main result,
Theorem~\ref{thm:main}.

\begin{theorem}
\label{thm:main2} \label{thm:equivalences}
Given a plane graph $G$, the following conditions are equivalent:
\begin{itemize}
\item[(i)] $G$ is generically rigid,
\item[(ii)] $G$ contains a spanning Laman subgraph,
\item[(iii)] $G$ can be labelled as a CPT with the generalized Laman property.
\item[(iv)] $G$ can be stretched as a pseudo-triangulation (with the given embedding and outer face).
\end{itemize}
\end{theorem}

The equivalence of (i) and (ii) is Laman's theorem and the fact that (i) and (iii) follow from (iv)
is Theorem \ref{thm:ordensantos}.
We will prove  (ii)$\Rightarrow$(iii)
(Section~\ref{sec:Laman-labelling}) and (iii)$\Rightarrow$(iv) (Section~\ref{sec:CPT-stretch}).

\begin{figure}[htb]
\centering
a. \includegraphics[height=1.3in]{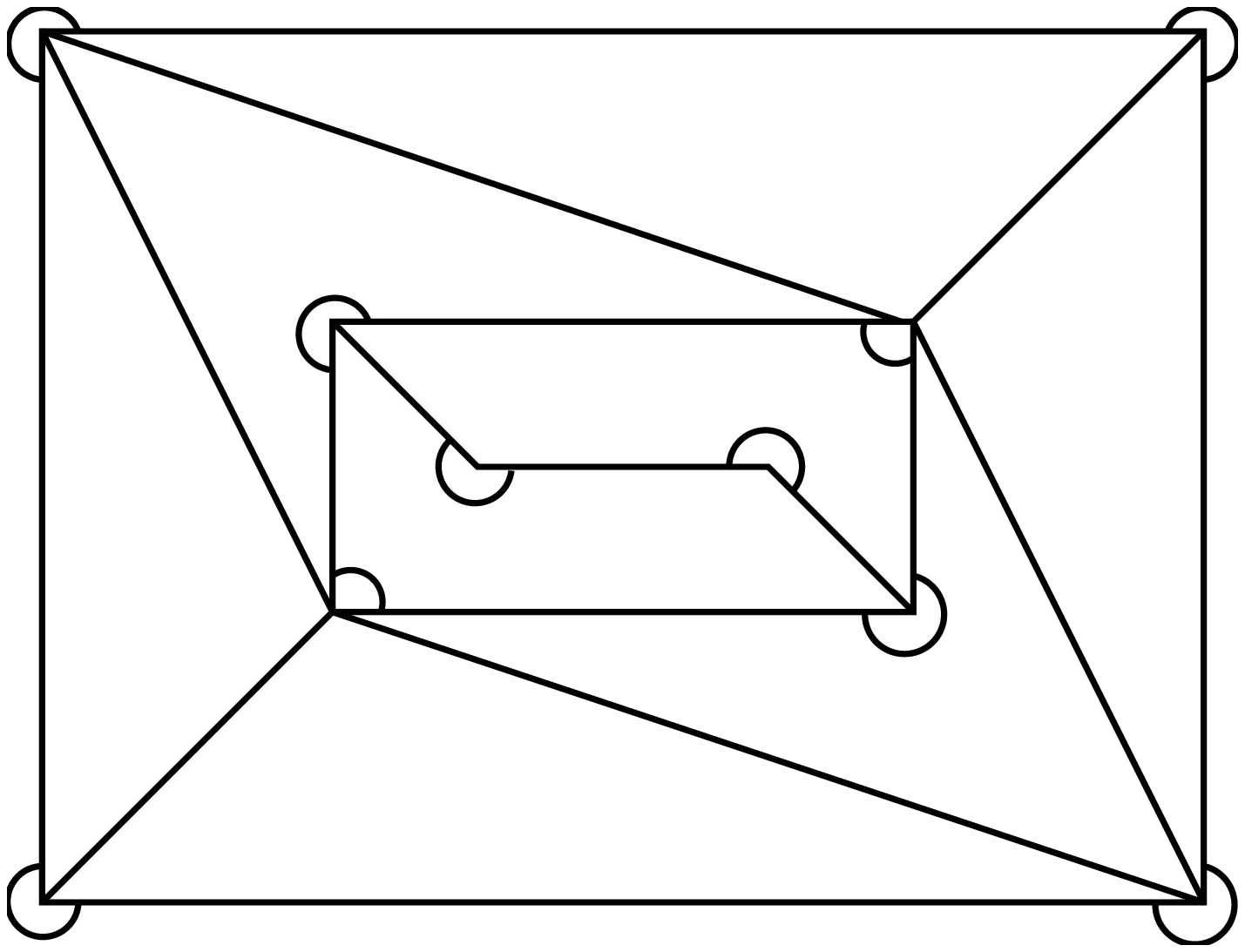}%,width=3.0in}
\qquad
b. \includegraphics[height=1.3in]{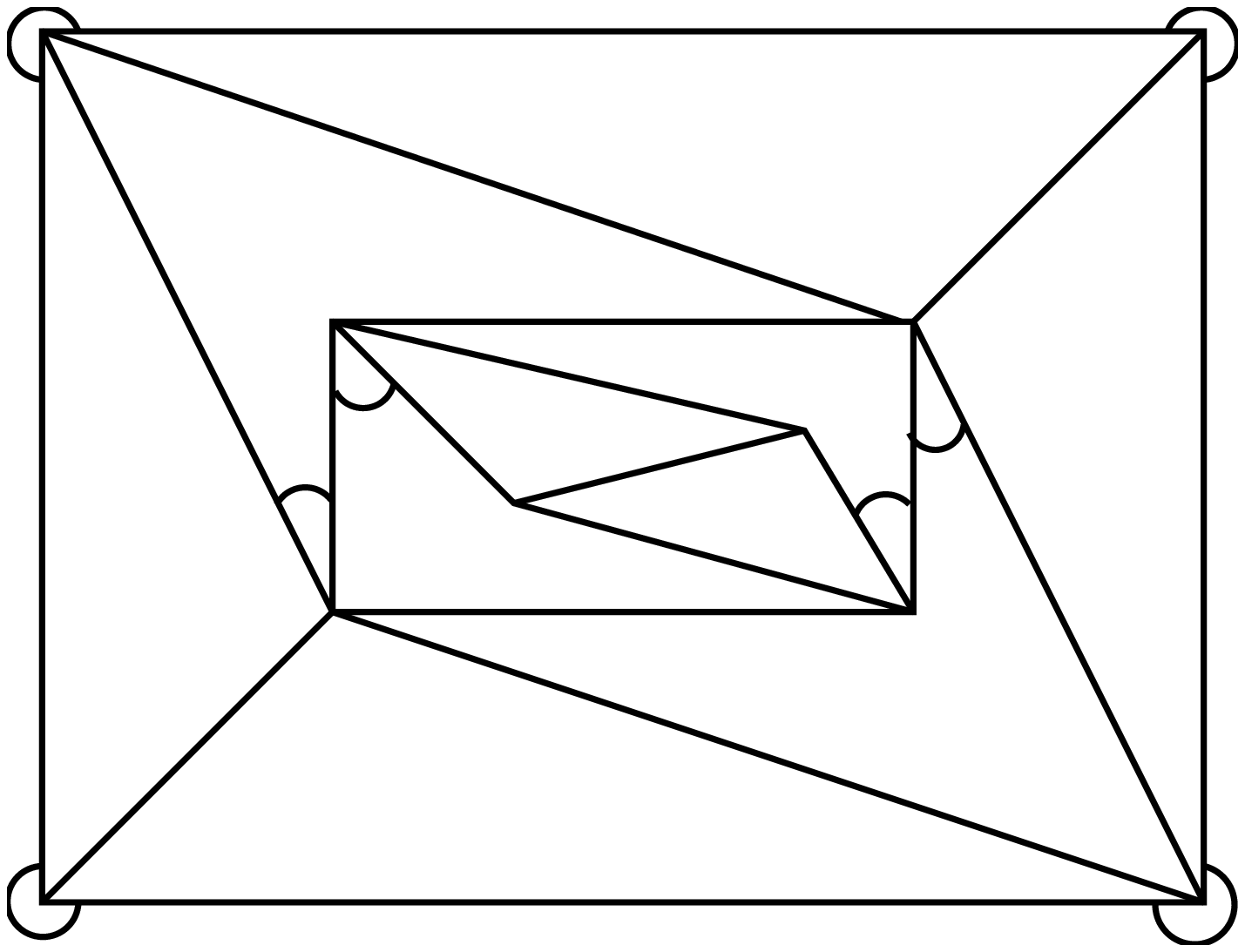}%,width=3.0in}
\caption{Left, a CPPT which is not a Laman graph. Right, a rigid
CPT which is not generalized Laman.\label{figure:NonLamanCPTs}}
\end{figure}

Note that having the generalized Laman property is not superfluous
in the statement, even for pointed CPT's.
Figure~\ref{figure:NonLamanCPTs}a shows a combinatorial pointed
pseudo-triangulation (CPPT)  which is not rigid because the innermost
three link chain has a motion or, equivalently, which is not Laman
because those four pointed vertices are incident to only seven
edges.

It is also not true that every rigid CPT has the generalized Laman
property, as Figure~\ref{figure:NonLamanCPTs}b shows, where the
two non-pointed vertices are incident to only five edges.

If we do not require the generalized Laman property, it is easy to
show that every rigid graph possesses a CPT labelling.  One can
start with a minimally rigid spanning subgraph, which has a CPPT
labelling by Theorem~\ref{thm:equivalences-lamangraphs}, and then
insert edges while only relabelling angles of the subdivided face.
For details see~\cite{orden}.  At each step one pointed vertex
must be sacrificed. But it is not obvious how to preserve the
generalized Laman condition in this process, even though one
starts out with a Laman graph.   In
Section~\ref{sec:Laman-labelling} we show that this can be done.

%%%%%%%%%%%%%%%%%%%%%%%%%%%%%%%%%%%%%%%%%%%%%%%%%%%%%%%%%%%%%%%
\section{Generalized Laman CPT's can be stretched}
\label{sec:CPT-stretch}

Here we prove the implication (iii)$\Rightarrow$(iv) of Theorem~\ref{thm:equivalences}. Our proof is
based on a partial result contained in Section~5 of \cite{horsssssw}. To state that result we need
to introduce the concept of corners of a subgraph.

Let $G=(V,E)$ be a CPT. Since $G$ comes (at least topologically) embedded in the plane,
we have an embedding of every subgraph of $G$. If $H$ is such a subgraph, every angle in $H$
is a union of one or more angles of $G$. Also, $H$ comes with a well-defined outer face, namely
the region containing the outer face of $G$. We say that a vertex $v$ of $H$ incident to the
outer face is a {\em corner} of $H$ if either
\begin{enumerate}
\item $v$ is pointed in $G$ and its big angle is contained in the outer face of $H$, or
\item $v$ is non-pointed in $G$ and it has two or more
consecutive small angles contained in the outer face of $H$.
\end{enumerate}

The following statement is Lemma~15 in \cite{horsssssw}:

\begin{lemma} \label{lemma:cornercount}
Let $H\subset G$ be a subgraph of a CPT and suppose that it is
connected and contains all the edges interior to its boundary
cycle (that is to say, $H$ is the graph of a simply connected
subcomplex of $G$).

Let $e$, $x$, $y$ and $b$ denote the numbers of edges, non-pointed
vertices, pointed vertices and length of the boundary cycle
in $S$, respectively. Then, the number $c_1$ of corners of the
first type (big angles in the outer boundary) of $H$ equals
\[
c_1 = e - 3x -2y + 3 + b.
\]
\end{lemma}

We say that a plane graph has \emph{non-degenerate faces} if the
edges incident to every face form a simple closed cycle. The
following statement is part of Theorem~7 of~\cite{horsssssw}:

\begin{theorem}
\label{thm:3cornerstretch}
For a combinatorial pseudo-triangulation $G$ with non-degenerate faces
the following properties are equivalent:
\begin{enumerate}
\item[(i)] $G$ can be stretched to become a pseudo-triangulation with the given assignment of angles.
\item[(ii)] Every subgraph of $G$ with at least three vertices has at least three corners.
\end{enumerate}
\end{theorem}

With this, in this section we only need to prove that:

\begin{theorem}
\label{thm:GLstretch} Let $G$ be a generalized Laman CPT.
Then:
\begin{enumerate}
\item[1.] Faces of $G$ are non-degenerate.
\item[2.] Every subgraph $H$  of $G$ on at least 3 vertices
has at least 3 corners.
\end{enumerate}
Hence, $G$ can be stretched.
\end{theorem}

\begin{proof}
1. Every face in a plane graph has a well-defined contour cycle.
What we need to prove is that no edge appears twice in the cycle.
For the outer face this is obvious, since all angles in it are
big: a repeated edge in the cycle would produce two big angles at
each of its end-points. Hence, assume that there is a repeated
edge $a$ in the contour cycle of a pseudo-triangle of $G$. This
implies that $G\setminus a$ has two components, ``one inside the
other''. Let us call $H$ the interior component. We will show that
the set of vertices of $H$ violates the generalized Laman property
by Lemma~\ref{lem:reformulated}.

Indeed, let $f$ and $e$ be the number of bounded faces and edges of $H$. Let
$x$ and $y$ be the numbers of non-pointed and pointed vertices in it.
The number of edges incident to the component is $e+1$ (for the edge $a$).
Hence, the generalized Laman property says that:
\[
e+1 \ge 3x + 2y.
\]
On the other hand, twice the number of edges of $H$ equals the number of angles of $G$
incident to $H$ minus one (because the removal of the edge $a$ merges two angles into one).
The number of small angles is at least $3f$ and the number of big angles is exactly $y$.
Hence,
\[
2e+1 \ge 3f+y.
\]
Adding these two equalities we get $3e+2 \ge 3f+3y+3x$, which violates
Euler's formula $e+1=f+y+x$.

2. Observe first that there is no loss of generality in assuming that
$H$ is connected (if it is not, the statement applies to each connected
component and the number of corners of $H$ is the sum of corners
of its components) and that $H$ contains all the edges of $G$  interior to
its boundary cycle (because these edges are irrelevant to the concept of corner).
We claim further that there is no loss of generality in assuming
that the boundary cycle of $H$ is non-degenerate. Indeed,
if $H$ has an edge $a$ that appears twice in its boundary cycle, its removal
creates two connected components $H_1$ and $H_2$, whose
numbers of vertices we denote $v_1$ and $v_2$. We claim that
each $H_i$ contributes at least $\min\{v_i,2\}$ corners to $H$.
Indeed, if $v_i$ is 1 or 2, then all vertices of $H_i$ are corners
in $H$. If $v_i\ge 3$, then $H_i$ has at least three corners and all
but perhaps one are corners in $H$.
Hence, $H$ has at least
$\min\{v_1,2\} + \min\{v_2,2\} \ge \min\{v_1+v_2,3\}$ corners, as desired.

Hence, we assume that $H$ consists of a simple closed cycle plus
all the edges of~$G$ interior to it.

Let $y$, $x$, $e$, and $b$ be the numbers of pointed vertices,
non-pointed vertices, edges and boundary vertices of $H$,
respectively. Let $V$ be the set of vertices of $H$ which are
either interior to $H$ or boundary vertices, but not corners. $V$
consists of $x+y-c_1-c_2$ vertices, where $c_1$ and $c_2$ are the
corners of type 1 and 2 of $H$, respectively. Hence,
Lemma~\ref{lem:reformulated}
 implies that the number of
edges incident to $V$ is at least
\[
2(y-c_1)+3(x-c_2) = 3x + 2y - 2c_1 - 3c_2.
\]
(Remark: $V$ certainly does not contain the corners of the whole
CPT, whose number is at least three.
This guarantees that we can apply Lemma~\ref{lem:reformulated} to
$V$).

On the other hand, the edges incident to $V$ are the $e-b$
interior edges of $H$ plus at most two edges per each boundary
non-corner vertex. Hence,
\[
3x + 2y - 2c_1 - 3c_2 \le e-b +2(b-c_1-c_2),
\]
or, equivalently,
\[
c_2 \ge 3x + 2y - e -b.
\]
Using Lemma \ref{lemma:cornercount} this gives $c_1+c_2\ge 3$.
\end{proof}

\begin{cor}
\label{cor:CPT-equivalences}
The following properties are equivalent for a combinatorial pseudo-triangulation $G$:
\begin{enumerate}
\item[(i)] $G$ can be stretched to become a pseudo-triangulation (with the given assignment of angles).
\item[(ii)] $G$ has the generalized Laman property.
\item[(iii)] $G$ has non-degenerate faces and
every subgraph of $G$ with at least three vertices has at least three corners.
\end{enumerate}
\end{cor}

\begin{proof}
The implications (i)$\Rightarrow$(ii),  (ii)$\Rightarrow$(iii),   and (iii)$\Rightarrow$(i),  are,
respectively, Theorems~\ref{thm:ordensantos}, \ref{thm:GLstretch} and \ref{thm:3cornerstretch}.
\end{proof}

To these three equivalences, Theorem~7 of \cite{horsssssw} adds a
fourth one: that a certain auxiliary graph constructed from $G$ is
3-connected in a directed sense. That property was actually the
key to the proof of (iii)$\Rightarrow$(i), in which the stretching
of $G$ is obtained using a directed version of Tutte's Theorem
saying that 3-connected planar graphs can be embedded with convex faces.

%%%%%%%%%%%%%%%%%%%%%%%%%%%%%%%%%%%%%%%%%%%%%%%%%%%%
\section{Obtaining Generalized Laman CPT-labellings}
\label{sec:Laman-labelling}

\begin{theorem}
The angles of a generically rigid plane graph can be labelled so that the
labelling is a CPT satisfying the generalized Laman condition.
\end{theorem}

The proof of the theorem proceeds by induction on the number $n$ of
vertices. As base case $n=3$ suffices. Consider a generically rigid
graph $G$ with more than three vertices. Since $G$ is generically
rigid it contains a Laman spanning subgraph~$L$. Let $v$ be a vertex
of minimal degree in $L$. Vertices in $L$ have minimum degree at
least two and $L$ has fewer than $2n$ edges (where $n$ is the number
of vertices), so $v$ has degree two or three:
\begin{itemize}
\item If $v$ has degree two in $L$, then $G\setminus v$ is generically rigid,
because $L\setminus v$ is a spanning Laman subgraph of it.
By inductive hypothesis, $G\setminus v$ has a generalized Laman
CPT labelling. Lemma~\ref{lem:degree2} shows that $G$ has one as well.
\item If $v$ has degree three in $L$, then $G\setminus v$ is either generically rigid
(and then we proceed as in the previous case) or it has one degree
of freedom. If the latter happens  there must be two neighbors $a$
and $b$ of $v$ such that if we insert the edge $e=ab$ into
$G\setminus v$ we get again a generically rigid graph. Since the plane
embedding of $G$ induces a plane embedding of $G':=G\setminus v
\cup e$, we have by inductive hypothesis a generalized Laman CPT
labelling of $G'$. Lemma~\ref{lem:degree3} shows that $a$ and $b$
can be chosen so that there is a CPT labelling of $G'$ extending
to one of $G$.
\end{itemize}

\begin{lemma}
\label{lem:degree2}
If $G\setminus v$ is generically rigid, then
every generalized Laman CPT labelling of
 $G\setminus v$ extends to one of $G$.
\end{lemma}

\begin{proof}
Let $T$ be the region of $G\setminus v$
where $v$ needs to be inserted. This region will either
be a pseudo-triangle  (of the CPT labelling of $G$) or
the exterior region of the embedding. For now, we assume that
$T$ is a pseudo-triangle. At the end of the proof we mention
how to proceed in the (easier) case where $T$ is the exterior
region.

Let $a$, $b$
and $c$ be the three corners of $T$. There are the following cases:

\begin{enumerate}
\item[(a)] If there is a pseudo-edge, say $ab$, containing at least two
 neighbors of $v$, then let $a'$ and $b'$ be the neighbors
closest to $a$ and $b$ on that pseudo-edge.
We may have $a=a'$ or $b=b'$, but certainly $a'\ne b'$.
We separate two subcases:
\begin{enumerate}
\item[(a.1)] If \emph{all} the neighbors of $v$ are on the pseudo-edge $ab$,
then let us label $v$ as pointed, with the big angle in the face containing $c$. The vertices
which are not neighbors of $v$ keep the angles they had in $G\setminus v$. The
neighbors of $v$ are labelled with both angles small, except that if $a\ne a'$
(respectively, $b\ne b'$)
then $a'$ (resp. $b'$) gets a big angle in the face containing $a$ (resp. $b$).
See Figure~\ref{fig:degree2}(a.1).

\item[(a.2)] If there is a neighbor of $v$ not on $ab$, let $c'$ be one of them,
chosen closest to $c$ (it may be $c$ itself). We
label $v$ as non-pointed, and the remaining angles as before except
if $c'\ne c$, we also put a big angle at $c'$, in the face containing $c$.
See Figure~\ref{fig:degree2}(a.2).
\end{enumerate}

\item[(b)] If no pseudo-edge contains two neighbors, then $v$ has either
two or three neighbors. The labelling is as in Figure~\ref{fig:degree2}(b),
with one of the edges from $v$ removed (and hence $v$ pointed) in the
case of only two neighbors. One of the  corners of $T$ may coincide
with the corresponding neighbor of $v$.
In this case that neighbor will get two small angles.

\begin{figure}[htb]
\centerline{
\includegraphics[scale=.75]{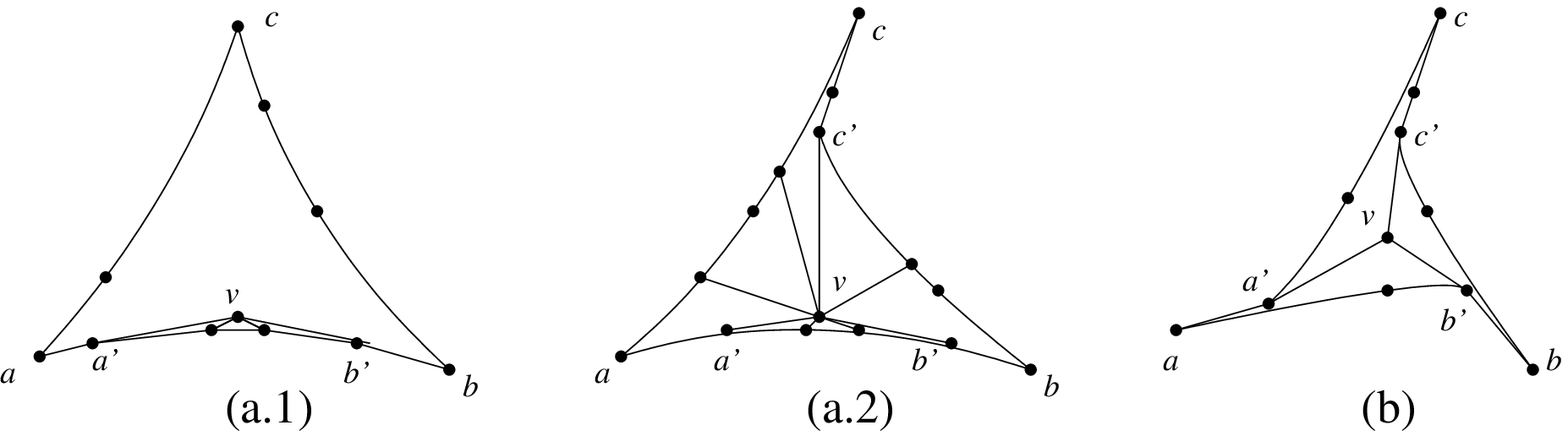}
}
\caption{}
\label{fig:degree2}
\end{figure}

\end{enumerate}

In all cases, the CPT labelling has the
property that every vertex that was non-pointed on $G\setminus v$
will remain non-pointed in $G$.
This in turn has the following
consequence: let $l$ be the number of neighbors of $v$
that were pointed in $G\setminus v$ and are not pointed in $G$.
By the count of edges in a CPT, the degree $d$ of
$v$ equals $l+2$ if $v$ is pointed and $l+3$ if $v$ is not pointed.
This relation follows also from the above case study. The
neighbors of $v$ that keep their status are precisely the points $a'$, $b'$ and
(if it exists) $c'$ in each case.

This is crucial in order to prove the generalized Laman property,
which we now do. Let $S$ be a subset of vertices of $G\setminus
v$. Since the subgraph induced by $S$ on $G$ and $G\setminus v$ is
the same, and since no vertex changed from non-pointed to pointed,
the generalized Laman property of $S$ in $G\setminus v$ implies
the same for $S$ in $G$. But we also need to check the property
for $S\cup v$. For this, by Lemma~\ref{lem:reformulated} it will
be enough if the number of neighbors of $v$ in $S$ that did not
change from pointed to non-pointed is at most two if $v$ is
pointed and at most three if $v$ is non-pointed. This follows from
the above equations $d=l+2$ and $d=l+3$ respectively.

As promised, we now address the case where $T$ is the exterior region of the
embedding of $G\setminus v$. This case can actually be considered a special
case of (a.1) above, since the exterior region has only ``one pseudo-edge''.
And, indeed, a labelling similar to the one shown in Figure~\ref{fig:degree2}(a.1)
works in this case, where the arc $a'b'$ now should be understood as the
segment of the boundary of the exterior region of $G\setminus v$ that becomes
interior in $G$.
\end{proof}

\begin{lemma}
\label{lem:degree3}
Let $v$ be a vertex of degree three in a Laman spanning subgraph of $G$
and let $e=ab$ be an edge between two neighbors of $v$ such that
$G':=G\setminus v \cup e$ is generically rigid.

\begin{enumerate}
\item Every CPT labelling of $G'$ extends to one of $G$.

\item If $a$ and $b$ are consecutive neighbors of $v$, then every generalized Laman
CPT labelling of $G'$ extends to a generalized Laman CPT labelling
of $G$.

\item If $a$ and $b$ are at ``distance two'' among neighbors of $v$, that is, if
there is a vertex $w$ such that $a$, $w$ and $b$
are consecutive neighbors of $v$, then either
\begin{enumerate}
\item[(i)] every generalized Laman
CPT labelling of $G'$ extends to a generalized Laman CPT labelling
of $G$, or
\item[(ii)] there is a generically rigid subgraph $H$ of $G\setminus v$ containing $w$
and with $v$ lying inside a bounded region of $H$.
\end{enumerate}
\end{enumerate}
\end{lemma}

It should be clarified what we mean by ``extends'' here. We mean
that all the labels of angles common to $G$ and $G'$ have the same
status. We exclude from this the angles at the end-points of the
edge $ab$, which may change status, even if these angles could in
principle be considered to survive in $G$, split into two edges
$av$ and~$vb$.

\begin{proof}

1. As in the previous Lemma, we start with a generalized Laman CPT
labelling of $G'=G\setminus v \cup e$. Let $T_1$ and $T_2$ be the
two pseudo-triangles containing~$e$. (As in the previous Lemma,
the case where one of them, say $T_1$, is the exterior region can
be treated as if $T_1$ was a pseudo-triangle with all neighbors of
$v$ in the same pseudo-edge). We call $a_i$, $b_i$ and $c_i$ the
vertices of $T_i$, in such a way that the pseudo-edge from $a_i$
to $b_i$ contains $a$ and $b$, in this order. Clearly, $a$
coincides with at least one of $a_1$ and $a_2$, and $b$ with one
of $b_1$ and $b_2$. Figure~\ref{fig:degree3-cases} shows the four
possibilities, modulo exchange of all $a$'s and $b$'s or all $1$'s
and $2$'s. In all cases we have drawn the union as a
pseudo-quadrilateral even if in parts (a) and (b) a
pseudo-triangle would in principle be possible. But that would
actually imply that $G\setminus v$ is generically rigid (because
it can be realized as a generalized Laman CPT), in which case we
could use Lemma~\ref{lem:degree2}.

\begin{figure}[htb]
\centerline{
\includegraphics[scale=0.75]{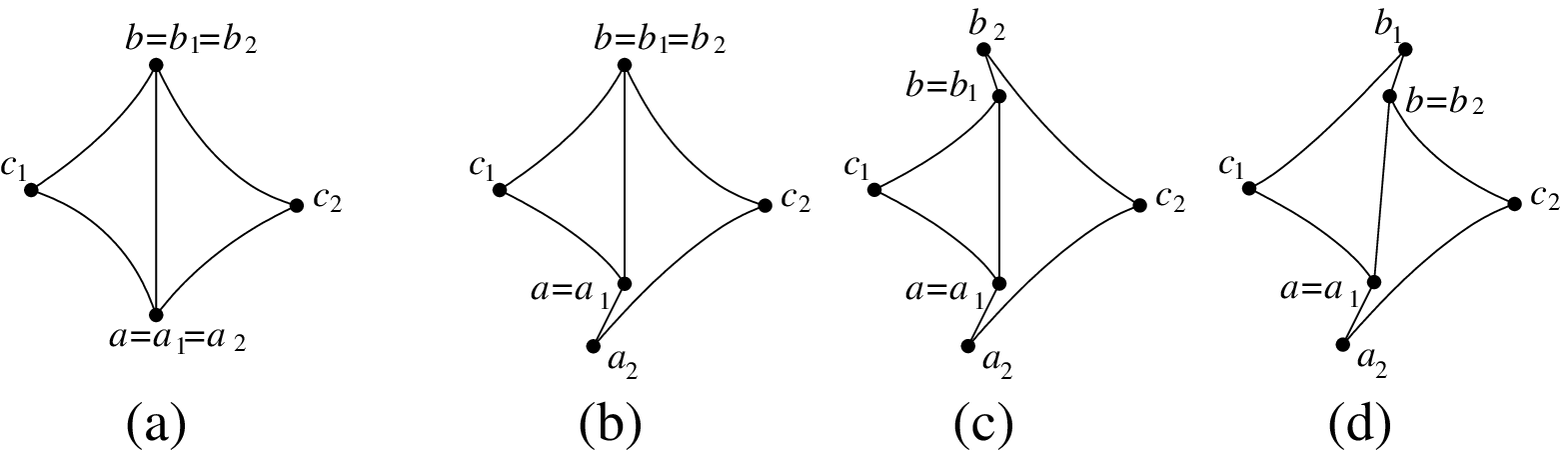}
}
\caption{}
\label{fig:degree3-cases}
\end{figure}

We will try to extend the labelling independently in $T_1$ and
$T_2$. We concentrate on one of them, say $T_1$. Let $a_1'$ and
$b'_1$ be the neighbors of $v$ closest to $a_1$ and $b_1$,
respectively, on the pseudo-edge $a_1b_1$. Observe that $a_1'$
will coincide with $a_1$ if $a_1$ is a neighbor, and will coincide
with $a$ if $a$  is the only neighbor on the path $aa_1$. If not
all the neighbors of $v$ in $T_1$ are on the pseudo-edge $a_1b_1$,
let $c_1'$ be one which is closest to $c_1$ (possibly $c_1$
itself). We assign labels as follows: All non-neighbors of $v$
keep their labels. All the angles at $v$ and at neighbors of $v$
are labelled small, with the following exceptions: If $a'_1\ne
a_1$, (respectively, $b'_1\ne b_1$ or $c'_1\ne c_1$) the angle at
$a'_1$ on the edge to $a_1$ is big (respectively, the angle at
$b'_1$ on the edge to $b_1$ or the angle at $c'_1$ on the edge to
$c_1$). Also, if all neighbors are on the pseudo-edge $a_1b_1$,
then the angle at $v$ on the pseudo-triangle $a_1b_1c_1$ is
labelled big. Figure~\ref{fig:degree3-labels} schematically shows
the two cases. For future reference, observe that the neighbors of
$v$ whose status does not change are precisely $a'_1$, $b'_1$ and,
unless $v$ gets a big angle,~$c'_1$.

\begin{figure}[htb]
\centerline{
\includegraphics[scale=.75]{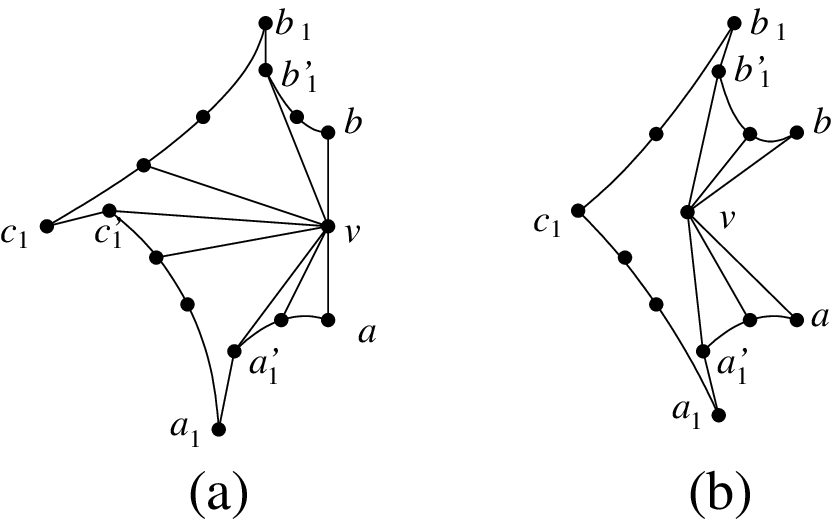}
}
\caption{}
\label{fig:degree3-labels}
\end{figure}

This clearly produces a pseudo-triangulation of $T_1$, and we use the same idea in $T_2$.
Observe also that we have not put big angles where there were none before.
In particular, no vertex other than perhaps $v$ will receive two big angles, and no vertex
that was not pointed in $G'$ will be pointed in $G$.

But $v$ itself may actually get two big angles, one on $T_1$ and one on $T_2$. This will happen if
all neighbors of $v$ in $T_1$ are on the pseudo-edge $a_1b_1$ and all neighbors in $T_2$
are on the pseudo-edge $a_2b_2$. In this case, since $v$ has at least three neighbors, one of $a'_1$,
$a'_2$, $b'_1$ or $b'_2$ is different from $a$ and $b$. Without loss of generality, suppose that $a'_2$ is different from
$a$. In particular $a_2\ne a$ and then $a_1=a=a'_1$ (as in parts (b), (c) or (d) of Figure \ref{fig:degree3-cases}).
In this case $a$ has received a small angle both in $T_1$ and in $T_2$, but it was originally pointed
in $G'$. We are then allowed to change the angle of $a$ in $T_1$ to be big, and that of $v$ in the same pseudo-triangle
to be small.
Figure~\ref{fig:degree3-change} shows the change. This finishes the proof of  part 1 of the Lemma.

\begin{figure}[htb]
\centerline{
\includegraphics[scale=.75]{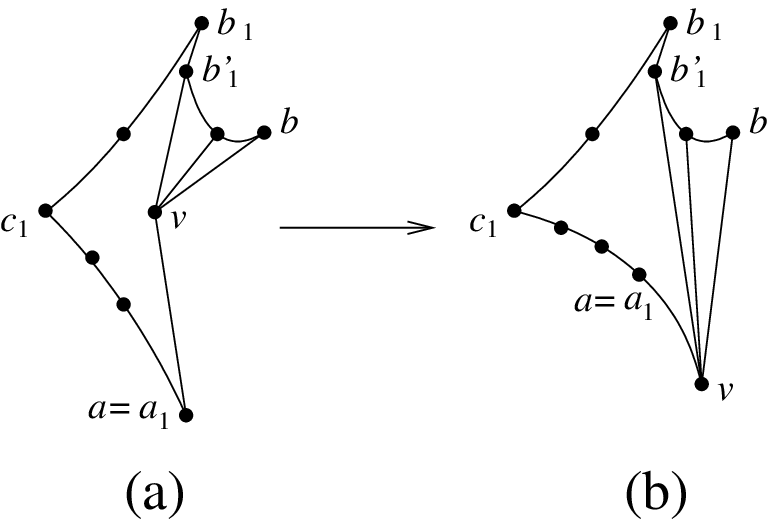}
}
\caption{}
\label{fig:degree3-change}
\end{figure}

Before going into the proofs of parts 2 and 3 we make the following observations about
the CPT labelling that we have just constructed:

\begin{enumerate}
\item[(a)] The number of neighbors of $v$ that do not change status is three if $v$ is pointed and
four if $v$ is non-pointed. This can be proved with a case study
using our explicit way of labelling, but it also follows from
global counts of small angles in the CPT's $G$ and $G'$.

\item[(b)] Assume now that $G'$ has the generalized Laman property, and let us try to
prove the property for $G$.
Every subset $S$ of vertices of $G$ not containing $v$ satisfies the generalized Laman count
in this CPT labelling: indeed, the subgraph induced by $G$
on $S$ is contained in the one induced by $G'$, and no vertex changed from non-pointed
to pointed. Hence, the Laman count translates from $G'$ to $G$.

\item[(c)] But if we try to prove the generalized Laman property for $S\cup v$, we encounter
a problem: Suppose that $a$ and $b$ do not both belong to $S$, so that the subgraph induced by $G$
on $S$ is the same as that induced by $G'$. Suppose also that $S$ is tight in $G'$,
meaning that the subgraph induced has exactly the number of edges permitted by the generalized Laman count.
Then, we need to prove that if $v$ is
pointed (respectively, non-pointed) at most two (respectively, three) of the neighbors of
$v$ in $S$ keep their pointedness status. But, globally, we know that three (respectively, four)
of the neighbors of $v$ keep their status, so we cannot finish the proof.

However, from this analysis we get very precise information on the
cases where $G$ happens not to have the generalized Laman
property. Namely, if $S\cup v$ fails to satisfy the generalized
Laman count, $S$ has the following four properties:
\begin{itemize}
\item $S$ does not contain both $a$ and $b$. Otherwise, the deleted edge $ab$ allows
for one extra edge to be inserted and the count is satisfied.
\item $S$ is tight in $G'$, otherwise again an extra edge is allowed to be inserted.
This  implies that if $G'$ is embedded as a pseudo-triangulation
(via Theorem~\ref{thm:GLstretch}), then the subgraph $G'|_S$
induced by $S$ is itself a pseudo-triangulation. Indeed, let
$e_S$, $x_S$ and $y_S$ be the numbers of edges, non-pointed
vertices and pointed vertices of $G'|_S$. Let $x$ and $y$ be the
numbers of vertices {\em of $S$} which are non-pointed and pointed
in $G'$, we have
\[
e_S=3x+2y-3\ge 3x_S+2y_S-3.
\]
Hence, $G'|_S$ is a pseudo-triangulation by
Lemma~\ref{lemma:pt-count}.

\item In
particular, $G'|_S$ is generically rigid. Since $G'|_S$ does not
contain $ab$, this implies that $S$ moves rigidly in the 1-degree
of freedom (\emph{1-dof}) mechanism $G\setminus v$.
\item $S$ contains all the neighbors of $v$ that did not change their status when we extended the CPT
labelling. In particular, at least one of $a$ or $b$ did change its status.
\end{itemize}
\end{enumerate}

2. We now suppose that $a$ and $b$ are consecutive neighbors of
$v$. Suppose that $T_1$ is the pseudo-triangle in which there is
no other neighbor. Certainly, $T_1$ does not impose any change of
status for $a$ or $b$. So, if $a$ or $b$ change their status, this
must be because of what happens in the pseudo-triangle $T_2$.
There are two cases:
\begin{itemize}
\item If our CPT-labelling of $T_2$ only changes the status of one of $a$ or $b$, we can restore its
status by the same type of change that we used when $v$ got two big angles:
in the pseudo-triangle $T_1$ we change the angle of $v$ from big to small, and that of the
neighbor that changed status in $T_2$ (which, clearly, had a big angle on $T_2$ and hence
a small angle in $T_1$), from small to big. We now have a CPT-labelling where neither $a$ nor
$b$ changed status, hence the generalized Laman property holds.
\item If our CPT-labelling of $T_2$ changes the status of both $a$ and $b$, let $a'$ and $b'$
be the neighbors of $v$ in the pseudo-edge $a_2b_2$ that do not
change their status by the CPT-labelling of $T_2$. We can try to
apply the trick of the previous case at vertex $a$ and at vertex
$b$. If the first one fails, we have a CPT-labelling where $a$, $a'$
and $b'$ did not change their status and without the generalized
Laman property. By our final remarks in the proof of part~1, this
implies that $a$, $a'$ and $b'$ move rigidly in the 1-dof mechanism
$G\setminus v$. Similarly, if the second fails, $b$, $b'$ and $a'$
move rigidly. Hence, if both fail, $a$, $b$, $a'$ and $b'$ move
rigidly in $G\setminus v$, which contradicts our initial choice of
the edge $e=ab$.
\end{itemize}

3. Applying the construction of part 1 to the edge $e$, it turns out
that $w$ must be one of the points that does not change status (that
is, one of the points $a'_i$, $b'_i$ or $c'_i$ of
Figure~\ref{fig:degree3-labels}). Suppose that the resulting CPT is
not generalized Laman, for some generalized Laman CPT of $G'$. We
want to prove that there is a rigid component of $G\setminus v$ that
includes $w$ and contains $v$ in the interior of a cycle.

By the remarks at the end of the proof of part 1, if the
generalized Laman property fails in $G$, then there is an induced
subgraph $S$ of $G\setminus v$ that contains all the neighbors of
$v$ that did not change status (in particular, contains $w$), and
which is itself a pseudo-triangulation (with respect to the CPT
labelling of $G'$).
 If this pseudo-triangulation already
contains $v$ in the interior of a pseudo-triangle, then the claim
is proved.
But it may
 happen that $v$ is in the exterior of this sub-pseudo-triangulation
 (see a schematic picture in Figure~\ref{fig:bad-pseudoquad}).

 \begin{figure}[htb]
\centerline{
\includegraphics[scale=.7]{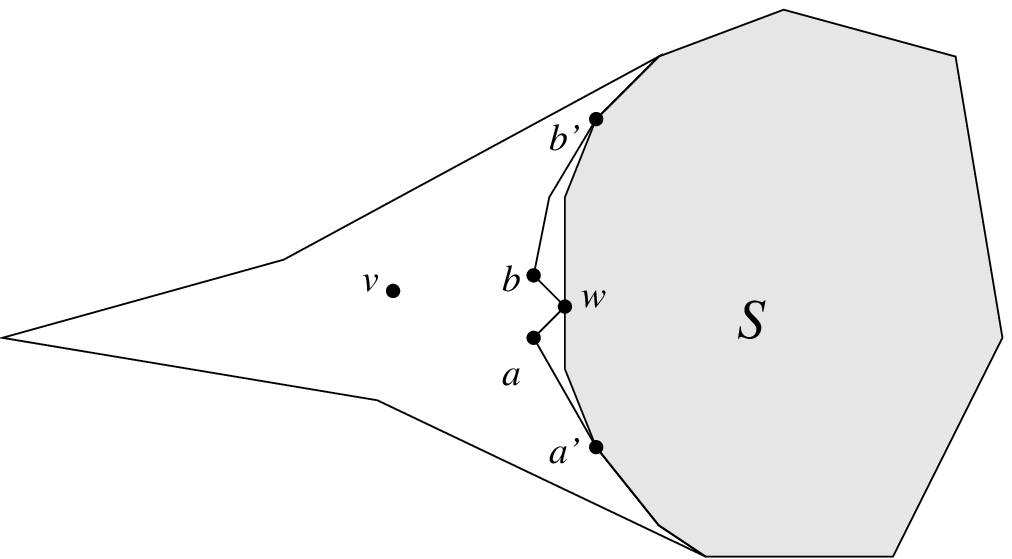}
}
\caption{}
\label{fig:bad-pseudoquad}
\end{figure}

 The crucial point now is that
 in the construction of part 1 the vertices that do not
 change status cannot all lie on the same
 pseudo-edge of the pseudo-quadrilateral of $G\setminus v$
containing $v$. Indeed, one
 of them is either $a'_1$ or $a'_2$
 (this is the point marked $a'$ in Figure~\ref{fig:bad-pseudoquad}),
 and lies on one of the two pseudo-edges ``on the $a$ side''. Another
 one is either $b'_1$ or $b'_2$
  (marked $b'$ in Figure~\ref{fig:bad-pseudoquad}),
 and lies on the two pseudo-edges ``on the $b$-side''.

 So, if $v$ is exterior to $S$, then the boundary of $S$ contains a
 concave chain connecting two different
 pseudo-edges of the pseudo-quadrilateral. Together with the opposite
 part of the pseudo-quadrilateral
 this produces a pseudo-triangle in $G\setminus v$ that contains $v$ in
 its interior and all the vertices
 that did not change status on its boundary. This pseudo-triangle,
 together with everything in its exterior,
 is a pseudo-triangulation, hence generically rigid
 (here, we are assuming that $G'$ and, in particular, $G\setminus v=G'\setminus e$
 has been stretched, via Theorem~\ref{thm:GLstretch}).
\end{proof}

\begin{cor}
\label{cor:degree3}
Let $v$ be a vertex of degree three in a Laman spanning subgraph of $G$. Then, there is an edge
$e=ab$ between two neighbors of $v$ such that
$G':=G\setminus v \cup e$ is generically rigid and  every generalized Laman
CPT labelling of $G'$ extends to a generalized Laman CPT labelling
of $G$.
\end{cor}

\begin{proof}
Let $w_1,\dots,w_k$ be the cyclic list of
neighbors of $v$ in $G$.
If there is an edge $e=w_iw_{i+1}$ (with indices meant modulo $k$)
between consecutive neighbors of
$v$ that restores rigidity in $G\setminus v \cup e$, then part 2
of the previous Lemma gives the statement for that $e$.
Hence, in the rest of the proof we suppose that this is not
the case.

We now look at edges of the form $e=w_{i-1}w_{i+1}$.
At least one of them must restore rigidity
in $G\setminus v \cup e$: if not, let
$G_0$ be the graph consisting of the $2k$ edges
$w_iw_{i+1}$ and $w_{i-1}w_{i+1}$. Our hypothesis is that
$G\setminus v \cup G_0$ is still a 1-dof mechanism. But, since
$G_0$ is rigid and contains all neighbors of $v$, this implies that
$G \cup G_0$, hence $G$, is a 1-dof mechanism too, a contradiction.

Let $S$ be the set of neighbors of $v$ in $G$ with the property that
$w_{i-1}w_{i+1}$ restores rigidity. In the rest of the proof we show
that there is a $w_i\in S$ such that $e=w_{i-1} w_{i+1}$  is as
claimed in the statement. We argue by contradiction, so assume the
claim is not true. Part 3 of the previous Lemma says that then for
each  $w_i\in S$ there is a generically rigid subgraph of
$G\setminus v$, let us denote it $H_{w_i}$,  such that $w_i\in
H_{w_i}$ and $v$ lies in the interior of a bounded face of
$H_{w_i}$. We now claim that the same holds for every $S'\subset S$:
there is a generically rigid subgraph $H_{S'}$ of $G\setminus v$
such that (1) $S'\subset H_{S'}$ and (2) $v$ lies in the interior of
a bounded face of $H_{S'}$. Indeed, after we know this for
one-element subsets we just need to show that from $H_{S_1}$ and
$H_{S_2}$ we can construct $H_{S_1\cup S_2}$. We consider first the
union of the two graphs $H_{S_1}$ and $H_{S_2}$. It clearly contains
both $S_1$ and $S_2$, and $v$ lies in the interior of a face: the
intersection of the faces $F_1$ of $H_{S_1}$ and $F_2$ of $H_{S_2}$
that contain $v$. The only problem is that $H_{S_1}\cup H_{S_2}$ may
not be generically rigid. Since $H_{S_1}$ and $H_{S_2}$ clearly
intersect (the boundaries $\partial F_1$ and $\partial F_2$ of $F_1$
and $F_2$ are two cycles around $v$ that must intersect because they
both contain neighbors of $v$ in $G$), if their union is not
generically rigid then they intersect in a single point. This point
must actually be in the two cycles
 $\partial F_1$ and $\partial F_2$. That is, the cycles
``are tangent and one is inside the other''. Since $S_1\subset \partial F_1$ and $S_2\subset \partial F_2$
consist only of neighbors of $v$ in $G$, this implies that one of $S_1$ and $S_2$ (say $S_1$) consists of a single
point $w_i$, which is the intersection point. In particular, $S_1\subset H_{S_2}$ and we can take
$H_{S_1\cup S_2}=H_{S_2}$.

So, taking $S'=S$, the conclusion is that all the $w_i$'s such that $w_{i-1}w_{i+1}$ restores
rigidity lie in a rigid subgraph of $G\setminus v$. Our final
claim is that under these conditions all neighbors of $v$ move
rigidly in the 1-dof mechanism  $G\setminus v$, which (as above)
contradicts the fact that $G$ is rigid. Indeed, assume without
loss of generality that the vertex $w_1$ is such that $w_kw_2$
restores rigidity. To seek a contradiction, let $w_{i+1}$ be the
first neighbor of $v$ (in the order $w_1,\dots,w_k$) which does
not move rigidly with $w_1$ and $w_2$, and let $w_{j+1}$ be the
first neighbor after $w_i$ which does not move rigidly with
$w_{i}w_{i+1}$. In particular, both $w_{i-1}w_{i+1}$ and
$w_{j-1}w_{j+1}$ restore rigidity and, by the above conclusion,
the three vertices $w_1$, $w_i$ and $w_j$ lie in a rigid subgraph.
This subgraph has two vertices in common with both
$\{w_1,\dots,w_i\}$ and $\{w_i,\dots,w_j\}$, which lie
respectively in two rigid subgraphs by the choice of $w_i$ and
$w_j$. Hence, all $\{w_1,\dots,w_j\}$ lies in a rigid subgraph, in
contradiction with the choice of $w_i$.
\end{proof}

%%%%%%%%%%%%%%%%%%%%%%%%%%%%%%%%%%%%%%%%%%%%%%%%%%%%%%%%%%%%%%%%
\section{Pseudo-triangulations on closed surfaces}
\label{sec:surfaces}

This section contains a couple of observations on the concept of
pseudo-triangula\-tions (and combinatorial ones) on closed
surfaces. We believe it would be interesting to develop this
concept further.

Let $G$ be a graph embedded on some closed surface of genus $g$.
  Every closed surface can be realized as the quotient of the sphere, the Euclidean
  plane, or the hyperbolic plane, by a discrete group of isometries, so in each case
  there is a well defined notion of distance and angle.
  Similar to the situation in the plane, a pseudo-triangulation of the surface is
a graph embedding with geodesic arcs such that every face has exactly three
angles smaller than $180^\circ$.  A combinatorial pseudo-triangulation is a
(topological) embedding together with an assignment of ``big'' and ``small'' to
angles such that every face has exactly three small angles and every vertex has
at most one big angle. One difference with the situation in the plane is that now
there is no ``outer'' face.

\begin{prop}\label{countProp}
    Let $G$ be embedded on a surface $S$ of genus $g$.
     If $G$ possesses a combinatorial
pseudo-triangular assignment then the numbers
 $e$, $x$ and $y$
    of edges, non-pointed vertices and pointed vertices satisfy
$e=3x +2y-6+6g$ if $S$ is
orientable, and $e=3x+2y-6+3g$ if $S$ is non-orientable.
\end{prop}

\begin{proof}
The number of angles equals twice the number of edges. There are $2e-y$ small
angles. The number of small angles equals three times the number of faces, $f$, which together with
Euler's formula $x+y-e+f=2-2g$ in the orientable case, or $x+y-e+f=2-g$ in the non-orientable case yields
the desired relationship between $x$, $y$ and $e$.
\end{proof}

\begin{figure}[htb]
\centering
\includegraphics[scale=.8]{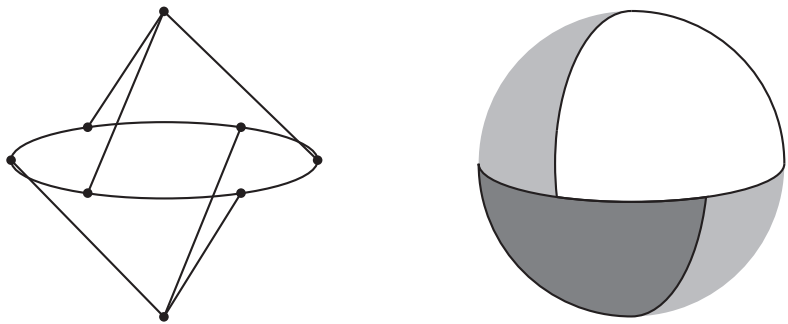}
\caption{A pseudo-triangular embedding of a cube on the sphere.\label{fig:pcube}}
\end{figure}

    The tree in Figure~\ref{FigMercedes}a is an example of a pointed pseudo-triangulation of the sphere.
    A triangular prism which can realized as a pointed pseudo-triangulation in the plane, cannot be realized as
    a pointed pseudo-triangulation of the sphere since, by Proposition~\ref{countProp}, any CPT labelling
    of the prism for the sphere must have three non-pointed vertices.  The graph of a cube, which has no
    CPT labelling for the plane, since it has too few edges to be rigid, can be realized as a pseudo-triangulation
    of the sphere, with two non-pointed vertices, see Figure~\ref{fig:pcube}
in which the pointed vertices are placed on the equator, and the two non-pointed vertices are at the poles.
In this geometric realization the large angles are exactly $180^\circ$,
so that pseudo-triangles are actual triangles,
although the complex is not a triangulation since it is not regular.

\begin{figure}[htb]
\centering
\includegraphics[scale=.85]{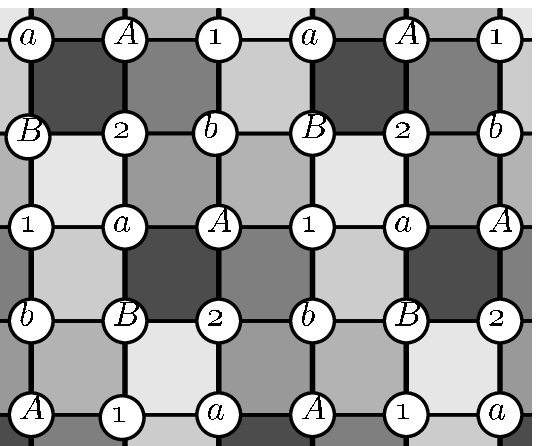}
\qquad
\includegraphics[scale=.85]{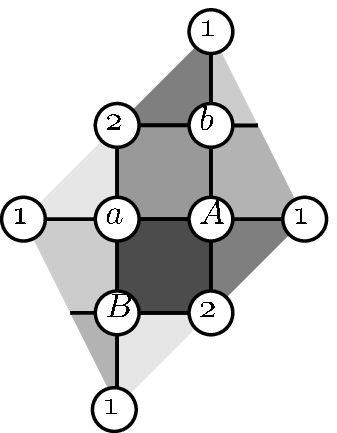}
\caption{The octahedron graph as a square tiling of the torus.\label{Figpseudo-tess1}}
\end{figure}

\begin{figure}[htb]
\centering
\includegraphics[scale=.85]{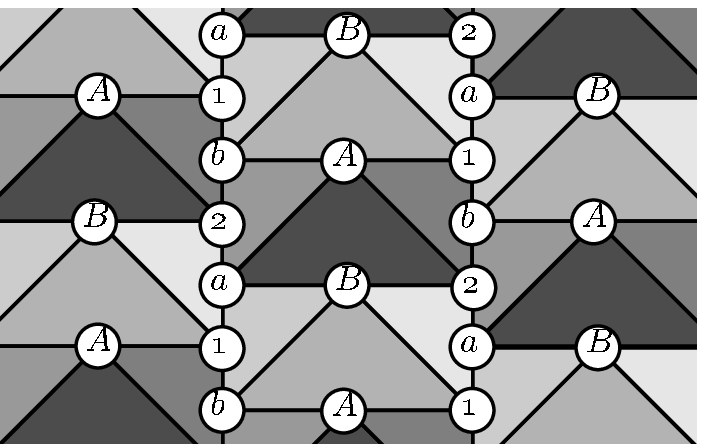}
\qquad
\includegraphics[scale=.85]{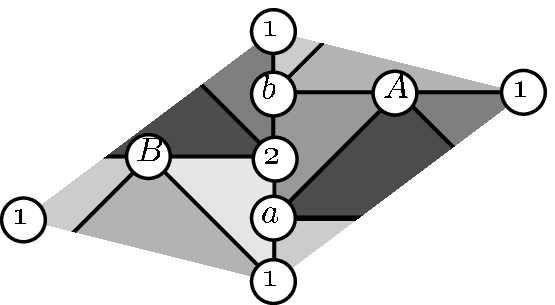}
\caption{The octahedron graph as a pointed pseudo-triangulation of the torus.\label{Figpseudo-tess2}}
\end{figure}

In Figure~\ref{Figpseudo-tess1} we see the well known embedding of the one-skeleton of the octahedron into
the torus with all square faces.
In Figure~\ref{Figpseudo-tess2}
we have modified this construction to embed the octahedron graph as a pointed pseudo-triangulation of the torus.
The octahedron graph, which is overbraced as a framework in the plane,
has no pointed pseudo-triangular embedding in either the plane or the sphere.

\section*{Acknowledgements}
This research was initiated at the Workshops on Topics in Computational Geometry organized
by Ileana Streinu at the Bellairs Research  Institute of McGill University in Barbados, Jan. 2002
and 2003, and partially supported by NSF Grant CCR-0203224. We also thank G.~Rote and an anonymous
referee for comments that helped us improve the
presentation.

%%%%%%%%%%%%%%%%%%%%%%%%%%%%%%%%%%%%%%%%%%%%%%%%%%%%%%%%%%%%%%%%%%%%%%%5
% Bibliography
%%%%%%%%%%%%%%%%%%%%%%%%%%%%%%%%%%%%%%%%%%%%%%%%%%%%%%%%%%%%%%%%%%%%%%%5


\begin{thebibliography}{11}

\bibitem{cegges}
{\sc B.\ Chazelle, H.\ Edelsbrunner, M.\ Grigni, L.\ Guibas, J.\ Erschberger, M.\ Sharir},
{\em Ray shooting in polygons using geodesic triangulations,}
Algorithmica {\bf  12},  no. 1, 54--68, 1994.


\bibitem{ChristopherServatius}
{\sc P.\ Christopher, B.\ Servatius,}
{\em Construction of self-dual graphs,}
American Mathematical Monthly, {\bf 99}, 2,  153--158, 1992.

\bibitem{gss}
{\sc J.\ Graver, B.\ Servatius and H.\ Servatius,}
``Combinatorial Rigidity'', AMS, Graduate Studies in Mathematics vol.2, 1993.

\bibitem{horsssssw}
{\sc R.\ Haas, D.\ Orden, G.\ Rote, F.\ Santos, B.\ Servatius, H.\ Servatius, D.\ Souvaine,
  I.\ Streinu, and W.\ Whiteley,}
{\em Planar minimally rigid graphs and pseudo-triangulations,}
Comput. Geom., 31, no.~1-2, 31--61, 2005.
% To
%appear in Computational Geometry, Theory and Applications, 2005.
%(Special issue for the Proc. 19th ACM Symp. Comp. Geom, 2003).

\bibitem{Laman}
{\sc G.\ Laman,}
{\em On graphs and rigidity of plane skeletal structures,}
J. Eng. Math., 4, 331–-340,
1970.

\bibitem{Maxwell}
{\sc J.\ C.\ Maxwell,}
{\em On reciprocal figures and diagrams of forces,}
Phil. Mag. Series 4, 27,  250--261, 1864.

\bibitem{orden}
{\sc D.\ Orden,}
{\em Two problems in geometric combinatorics: Efficient triangulations of the hypercube; Planar graphs and rigidity,}
Ph.D. Thesis, Universidad de Cantabria, 2003.

\bibitem{osssw}
{\sc D.\ Orden, G.\ Rote, F.\ Santos, B.\ Servatius, H.\
Servatius, W.\ Whiteley} {\em Non-crossing frameworks with
non-crossing reciprocals,} Discrete and Computational Geometry
(special issue in honor of Lou Billera), {\bf 32}:4 (2004),
567--600.

\bibitem{os-pncgp-02}
{\sc D.~Orden and F.~Santos,} {\em The polytope of non-crossing
graphs on a planar point set}, Discrete and Computational Geometry
{\bf 33}:2 (2005), 275--305.
\bibitem{pv}
{\sc M.\ Pocchiola and G.\ Vegter,}
{\em Topologically sweeping visibility complexes via pseudo-triangulations,}
Discrete and Computational Geometry, 16, 419--453, 1996.


\bibitem{rss}
{\sc G.\ Rote, F.\ Santos, and I.\ Streinu}
{\em Expansive motions and the polytope of pointed pseudo-triangulations,}
in ``Discrete and Computational Geometry -- The Goodman-Pollack
     Festschrift'', (B.\ Aronov, S.\ Basu, J.\ Pach, M.\ Sharir, eds.),
Algorithms and Combinatorics vol. 25, Springer Verlag, Berlin,  699--736, 2003.
%({\tt http://arXiv.org/abs/math.CO/0206027})

\bibitem{streinu}
{\sc I.\ Streinu,}
{\em A combinatorial approach to planar non-colliding robot arm motion
     planning,}
In ``Proc. 41st Ann. Symp. on Found. of Computer Science'' (FOCS
     2000), Redondo Beach, CA, 443--453, 2000.


\bibitem{walterS}
{\sc W.\ Whiteley,}
{\em Rigidity and Scene Analysis},
in ``Handbook of Discrete and Computational Geometry'',
(J.\ E.\ Goodman, J.\ O'Rourke, eds.), 893--916, 1997.

\end{thebibliography}
\end{document}